\documentclass[11pt]{amsart}
\usepackage{amsmath,amssymb}

\usepackage{mathrsfs,bbm,eucal}

\theoremstyle{plain}
\newtheorem{thm}{Theorem}
\theoremstyle{plain}
\newtheorem{lem}{Lemma}
\theoremstyle{plain}
\newtheorem{prop}{Proposition}
\theoremstyle{definition}
\newtheorem{defn}{Definition}
\theoremstyle{remark}

\newcommand{\dirac}{\mbox{$\mathcal{D}\!\!\!\!\!\:/\!\;$}}

\newcommand{\spinor}{\mbox{$S\!\!\!\!\!\:/\;\!$}}
\newcommand{\bundle}[1]{\CMcal{#1}}

\newcommand{\R}{\mathbbm{R}}
\newcommand{\C}{\mathbbm{C}}

\newcommand{\T}{\CMcal{T}}

\newcommand{\hyp}[1]{\mbox{$#1\mathrm{H}$}}

\begin{document}
\title[Scalar curvature rigidity of $\hyp{\C }^{2n+1}$]{Scalar curvature rigidity of almost Hermitian spin manifolds which are asymptotically complex hyperbolic}
\author{Mario Listing}
\address{Department of Mathematics, Stony Brook University, Stony Brook, NY 11794-3651, USA}
\email{ listing@math.sunysb.edu}
\thanks{Supported by the German Research Foundation}
\begin{abstract}
This paper generalizes a rigidity result of complex hyperbolic spaces by M.~Herzlich. We prove that an almost Hermitian spin manifold $(M,g)$ of real dimension $4n+2$ which is strongly asymptotic to $\hyp{\C }^{2n+1}$ and satisfies a certain scalar curvature bound must be isometric to the complex hyperbolic space. The fact that we do not assume $g$ to be K\"ahler reflects in the inequality for the scalar curvature.
\end{abstract}
\keywords{almost complex structures, rigidity, K\"ahler Killing spinors}
\subjclass[2000]{Primary 53C24, Secondary 53C55}
\maketitle

\section{Introduction}
Rigidity of symmetric spaces of non--compact type is a frequently studied problem (cf.~\cite{Bart,He3,MinO,MinO2}). Based on E.~Witten's idea in the proof of the positive energy theorem (cf.~\cite{Wit}), R.~Bartnik showed in \cite{Bart} that an asymptotically flat spin manifold of non--negative scalar curvature and with vanishing mass must be the Euclidean space. The analogous rigidity result for the real hyperbolic was proved by M.~Min--Oo in \cite{MinO}, in particular a strongly asymptotically hyperbolic spin manifold $(M^n,g)$ with scalar curvature $\mathrm{scal}\geq -n(n-1)$ is isometric to the hyperbolic space. Moreover, M.~Herzlich showed in \cite{He3} that a strongly asymptotically complex hyperbolic K\"ahler spin manifold $(M^{2m},g)$ of odd complex dimension $m$ and with scalar curvature $\mathrm{scal}\geq -4m(m+1)$ must be isometric to the complex hyperbolic space $\hyp{\C}^m$.

In this paper we generalize Herzlich's result in the way that we replace the K\"ahler assumption by the weaker condition: almost Hermitian. 

\begin{defn}
\rm $(\hyp{\C}^m,g_0)$ denotes the complex hyperbolic space of complex dimension $m$ and holomorphic sectional curvature $-4$, i.e.~$K\in [-4,-1]$, as well as $B_R(q)\subset M$ is the set of all $p\in M$ with geodesic distance to $q$ less than $R$. Let $(M^{2m},g,J)$ be an almost Hermitian manifold, i.e.~$g$ is a Riemannian metric and $J$ is a $g$--compatible almost complex structure. $(M,g,J)$ is said to be \emph{strongly asymptotically complex hyperbolic} if there is a compact manifold $C\subset M$ and a diffeomorphism $f:E:=M-C\to \hyp{\C}^m-\overline{B_R(0)}$ in such a way that the positive definite gauge transformation $A\in \Gamma (\mathrm{End}(TM_{|E}))$ given by
\[
g(AX,AY) = (f^{*}g_{0})(X,Y)\quad g(AX,Y) = g(X,AY)
\]
satisfies:
\begin{enumerate}
\item $A$ is uniformly bounded.
\item Suppose $r$ is the $f^*g_0$--distance to a fixed point, $\nabla ^0$ is the Levi--Civita connection for $f^*g_0$ and $J_0$ is the complex structure of $\hyp{\C}^m$ pulled back to $E$, then \[
\left| \nabla ^{0}A\right| +\left| A-\mathrm{Id} \right| +\left| AJ_0-J\right| \in L^{1}(E;e^{2r}\mathrm{vol}_{g})\cap L^{2}(E;e^{2r}\mathrm{vol}_{g}).
\]
\end{enumerate}
\end{defn}
In particular, in contrast to the previous definition and result by Herzlich, a compact conformal transformation of the standard metric on $\hyp{\C}^m$ supplies a manifold which is strongly asymptotically complex hyperbolic.
\begin{thm}
\label{rig_thm}
Let $(M^{4n+2},g,J)$ be a complete almost Hermitian spin manifold of odd complex dimension $m=2n+1$. If $(M,g,J)$ is strongly asymptotically complex hyperbolic and satisfies the scalar curvature bound
\begin{equation}
\label{ineq}
\mathrm{scal}\geq -4m(m+1)+2\Bigl[ |\mathrm{d}^*\Omega |+|\CMcal{D}^\prime \Omega |+|\CMcal{D}^{\prime \prime }\Omega |\Bigl] ,
\end{equation}
then $(M,g,J)$ is K\"ahler and isometric to $\hyp{\C}^m$. 
\end{thm}
In this case $\Omega =g(.,J.)$ is the $2$--form associated to $J$, $\mathrm{d}^*$ is formal $L^2$--adjoint of the exterior derivative $\mathrm{d}$ and $\CMcal{D}^\prime +\CMcal{D}^{\prime \prime }$ is the Dolbeault decomposition of $ \CMcal{D}=\mathrm{d}+\mathrm{d}^*$  in $\Lambda ^*(TM)\otimes \C $, i.e.~if $e_1,\ldots ,e_{2m}$ is an orthonormal base, we define $\CMcal{D}^\prime =\sum e_j^{1,0}\cdot \nabla _{e_j}$ and $\CMcal{D}^{\prime \prime }=\sum e_j^{0,1}\cdot \nabla _{e_j}$. Introduce $\CMcal{D}^c:=\mathrm{d}^c+\mathrm{d}^{c,*}$ with $\mathrm{d}^c:=\sum J(e_k)\wedge \nabla _{e_k}$ and $\mathrm{d}^{c,*}:=-\sum J(e_k)\llcorner \nabla _{e_j}$, we obtain $\CMcal{D}^\prime =\frac{1}{2}(\CMcal{D}-\mathbf{i}\CMcal{D}^c)$ as well as $\CMcal{D}^{\prime \prime }=\frac{1}{2}(\CMcal{D}+\mathbf{i}\CMcal{D}^c)$. In particular, we can estimate
\[
|\CMcal{D}^\prime \Omega |+|\CMcal{D}^{\prime \prime }\Omega |\leq |\mathrm{d}^*\Omega |+|\mathrm{d}\Omega |+|\mathrm{d}^{c,*}\Omega |+|\mathrm{d}^c\Omega |.
\]

The proof of this rigidity theorem is as usual based on the non--compact Bochner technique which was introduced by Witten in \cite{Wit}. We show an integrated Bochner--Weitzenb\"ock formula for the K\"ahler Killing connection which allows the usage of this technique. We expect to prove a similar result in the complex even--dimensional case and for the quaternionic hyperbolic space, but because of representation theoretical problems, there will be more terms involved in inequality (\ref{ineq}).
\section{Preliminaries}
Let $(M,g,J)$ be an almost Hermitian spin manifold of complex dimension $m$ and denote by $\gamma $ respectively $\cdot $ the Clifford multiplication on the complex spinor bundle $\spinor M$ of $M$. $\spinor M$ decomposes orthogonal into
\begin{equation}
\label{dec_spinor}
\spinor M=\spinor _0\oplus \cdots \oplus \spinor _m
\end{equation}
(cf.~\cite{Kir1,LaMi}) where each $\spinor _j$ is an eigenspace of $\Omega =g(.,J.)$ to the eigenvalue $\mathbf{i}(m-2j)$. We denote by $\pi _j$ the orthogonal projection $\spinor M\to \spinor _j$. The decomposition (\ref{dec_spinor}) is parallel (i.e.~$\nabla \pi _j=0 $ for all $j$) if $(g,J)$ is K\"ahler. As usual we introduce $X^{1,0}:=\frac{1}{2}(X-\mathbf{i}J(X))$ as well as $X^{0,1}:=\frac{1}{2}(X+\mathbf{i}J(X))$ and obtain $\gamma (X^{1,0}):\spinor _j\to \spinor _{j+1}$ as well as $\gamma (X^{0,1}):\spinor _j\to \spinor _{j-1}$, where $\spinor _j=\{ 0\} $ if $j\notin \{ 0,\ldots ,m\}$.

Supposing $(g,J)$ to be K\"ahler and $m=2n+1$ to be odd, then a K\"ahler Killing spinor (cf.~\cite{Kir1}) is a section in $\spinor _n\oplus \spinor _{n+1}$ which is parallel w.r.t.
\[
\nabla _X+\kappa \left( \gamma (X^{1,0})\pi _n+\gamma (X^{0,1})\pi _{n+1}\right) .
\]
In particular, if there is a non--trivial K\"ahler Killing spinor, $g$ is Einstein of scalar curvature $4m(m+1)\kappa ^2$. Moreover, the subbundle $\spinor _n\oplus \spinor _{n+1}$ is trivialized by K\"ahler Killing spinors on $\hyp{\C}^m$ if we choose $\kappa =\pm \mathbf{i}$.

\section{Bochner--Weitzenb\"ock formula}
Suppose $(M,g,J)$ is spin and almost Hermitian of odd complex dimension $m=2n+1$. We define $\bundle{V}:=\spinor _n\oplus \spinor _{n+1}$, its projection $\mathrm{pr}_\bundle{V}:=\pi _n+\pi _{n+1}$ and
\[
\frak{T}_X:=\mathbf{i}\left( \gamma (X^{1,0})\pi _{n}+\gamma (X^{0,1})\pi _{n+1}\right) .
\]
Since $(\gamma (X^{1,0})\pi _j)^*=-\gamma (X^{0,1})\pi _{j+1}$, $\frak{T}$ is a selfadjoint endomorphism on $\bundle{V}$ (respectively $\spinor M$). Define the connection $\widehat{\nabla }:=\nabla +\frak{T}$ on $\spinor M$. The Dirac operator of $\widehat{\nabla }$ is given by $\widehat{\dirac}=\dirac +\T $ where $\dirac $ is the Dirac operator of $\nabla $ and $\T$ equals
\[
-\mathbf{i}(m+1)\mathrm{pr}_{\bundle{V}}
\]
in this case we used (cf.~\cite{Kir1})
\begin{equation}
\label{rel1}
\sum _k e_k\cdot e_k^{1,0}=-m+\mathbf{i}\gamma (\Omega )\ \ \text{and}\ \ \sum _k e_k\cdot e_k^{0,1}=-m-\mathbf{i}\gamma (\Omega ).
\end{equation}
Since $\gamma (X)\T $ is not selfadjoint on the full spinor bundle, we consider instead $\mathbb{T}:=-\mathbf{i}(m+1)$ as well as the Dirac operator $\widetilde{\dirac}:=\dirac +\mathbb{T}$.
\begin{prop}
Let $(M,g,J)$ be almost Hermitian of odd complex dimension $m$, then the integrated Bochner--Weitzenb\"ock formula
\[
\int\limits _{\partial N}\left< \widehat{\nabla }_\nu \varphi +\nu \cdot \widetilde{\dirac}\varphi ,\psi \right> = \int\limits _N \left< \widehat{\nabla}\varphi ,\widehat{\nabla }\psi \right> -\left< \widetilde{\dirac}\varphi ,\widetilde{\dirac}\psi \right>+\left< \widehat{\frak{R}}\varphi ,\psi \right>
\]
holds for any compact $N\subset M$ and $\varphi ,\psi \in \Gamma (\spinor M)$. In this case $\nu $ is the outward normal vector field on $\partial N$ and $\widehat{\frak{R}}$ is given by
\[
\frac{\mathrm{scal}}{4}+m(m+1)+(m+1)^2\mathrm{pr}_{\bundle{V}^\perp}+\delta \frak{T} 
\]
while $\mathrm{pr}_{\bundle{V}^\perp}$ is the projection to the orthogonal complement of $\bundle{V}$ in $\spinor M $ and $\delta \frak{T}$ is the divergence of $\frak{T}$, i.e.~$\delta \frak{T}=\sum (\nabla _{e_j}\frak{T})_{e_j}$. Moreover, the boundary operator $\widehat{\nabla }_\nu +\nu \cdot \widetilde{\dirac}$ is selfadjoint.
\end{prop}
\begin{proof}
The essential facts are $(\frak{T}_X)^*=\frak{T}_X$ and $(\gamma (X)\mathbb{T})^*=\gamma (X)\mathbb{T}$. In particular, since $\nabla _\nu +\nu \cdot \dirac $ is a selfadjoint boundary operator,  $\widehat{\nabla }_\nu +\nu \cdot \widetilde{\dirac}$ is selfadjoint. The formal $L^2$--adjoint of $\widetilde{\dirac}$ is given by $\widetilde{\dirac}^*=\dirac -\mathbb{T}$. Thus, we can easily verify
\[
\int\limits _N \left< \widetilde{\dirac}\varphi ,\widetilde{\dirac}\psi \right> =-\int\limits _{\partial N}\left< \nu \cdot \widetilde{\dirac}\varphi ,\psi \right> +\int\limits _N\left< \widetilde{\dirac}^*\widetilde{\dirac}\varphi ,\psi \right> 
\]
as well as $\widetilde{\dirac}^*\widetilde{\dirac}=\dirac ^2+(m+1)^2$. Moreover, using $(\frak{T}_X)^*=\frak{T}_X$ on $\spinor M$ leads to
\[ \begin{split}
\int\limits _N\left< \widehat{\nabla }\varphi ,\widehat{\nabla }\psi \right> =&\int\limits _N \left< \nabla \varphi ,\nabla \psi \right> +\left< \nabla \varphi ,\frak{T}\psi \right> +\left< \frak{T}\varphi ,\nabla \psi \right> +\left< \frak{T}\varphi ,\frak{T}\psi \right>\\
=& \int\limits _{\partial N}\left< \nabla _\nu \varphi +\frak{T}_\nu \varphi ,\psi \right> +\int\limits _N \left< \nabla ^*\nabla \varphi ,\psi \right> +\\
&\quad +\int\limits _N\left< \frak{T}\varphi ,\frak{T}\psi \right> -\left< \delta \frak{T}\varphi ,\psi \right> 
\end{split}\]
for all $\varphi ,\psi \in \Gamma (\spinor M)$. We use the facts $\pi _j\gamma (X)\pi _{j-1}=\gamma (X^{1,0})\pi _{j-1}$ and $\pi _j\gamma (X)\pi _{j+1}=\gamma (X^{0,1})\pi _{j+1}$ as well as (\ref{rel1}) to compute
\[ \begin{split}
\left< \frak{T}\varphi ,\frak{T}\psi \right> =&\sum _k\left< e_k\cdot \varphi _n,e_k^{1,0} \cdot \psi _n\right> +\sum _k\left< e_k\cdot \varphi _{n+1},e_k^{0,1} \cdot \psi _{n+1}\right> \\
=&(m+1)\left< \mathrm{pr}_\bundle{V}\varphi ,\psi \right> .
\end{split}\]
In particular, the Lichnerowicz formula $\dirac ^2=\nabla ^*\nabla +\frac{\mathrm{scal}}{4}$ gives the claim with
\[
\widehat{\frak{R}}=\frac{\mathrm{scal}}{4}+(m+1)^2-(m+1)\mathrm{pr}_\bundle{V}+\delta \frak{T}.
\] 
\end{proof}
\begin{lem}
Suppose inequality (\ref{ineq}) of the main theorem holds, then at each point of $M$, $\widehat{\frak{R}}$ has no negative eigenvalues: $\widehat{\frak{R}}\geq 0$.
\end{lem}
\begin{proof}
We have to find an estimate for $\delta \frak{T}$. Let $e_1,\ldots ,e_{2m}$ be normal coordinates at $T_pM$ with $e_{m+j}:=Je_j$ in $p$. We obtain
\[
\begin{split}
\delta \frak{T}=&\sum _{j=1}^{2m}(\nabla _{e_j}\frak{T})_{e_j}\\
=&\frac{1}{2}\gamma (\delta J)(\pi _{n}-\pi _{n+1})+\mathbf{i}\sum _{j=1}^{2m}\left( e_j^{1,0}\cdot \nabla _{e_j}\pi _n+e_j^{0,1}\cdot \nabla _{e_j}\pi _{n+1}\right) .
\end{split}\]
Thus, we have to estimate $\nabla _X\pi _r$ for $r=n,n+1$. We conclude from $\pi _n\gamma (\Omega )=\mathbf{i}\pi _n$
\[
(\nabla _X\pi _n)(\mathbf{i}-\gamma (\Omega ))=\pi _n\gamma (\nabla _X\Omega )
\]
as well as from $\gamma (\Omega )\pi _n=\mathbf{i}\pi _n$
\[
(\mathbf{i}-\gamma (\Omega )) (\nabla _X\pi _n)=\gamma (\nabla _X\Omega )\pi _n .
\]
Using the facts $\pi _n (\nabla _X\pi _n)\pi _n=0$ and
$
\mathbf{i}-\gamma (\Omega )=\sum _{j\neq n}c_j\pi _j$
with $|c_j|\geq 2$, $|\nabla _X\pi _n|$ can be estimated by $\frac{1}{2}|\nabla _X\Omega |$. Thus,
\[
\sum _{j=1}^{2m}e_j^{1,0}\cdot (\nabla _{e_j}\pi _n)(\mathbf{i}-\gamma (\Omega ))=\pi _{n+1}\sum _{j=1}^{2m}\gamma (e_j^{1,0}\cdot \nabla _{e_j}\Omega )
\]
leads to
\[
\Biggl| \sum _{j=1}^{2m}e_j^{1,0}(\nabla _{e_j}\pi _n)\phi \Biggl|\leq \frac{1}{2}\bigl| \gamma (\CMcal{D}^\prime \Omega )\phi \bigl| ,
\]
if $\pi _n(\phi )=0$. Moreover,
\begin{eqnarray*}
\sum _{j=1}^{2m}\gamma (e_j^{1,0}\cdot \nabla _{e_j}\Omega )\pi _n&=& \sum _{j=1}^{2m}\gamma (e_j^{1,0})(\mathbf{i}-\gamma (\Omega ))(\nabla _{e_j}\pi _n)\\
&=& -\sum _{j=1}^{2m} (\mathbf{i}+\gamma (\Omega ))\gamma (e_j^{1,0}) (\nabla _{e_j}\pi _n)
\end{eqnarray*}
shows
\[
\Biggl| \sum _{j=1}^{2m}e_j^{1,0}(\nabla _{e_j}\pi _n)\phi \Biggl|\leq \frac{1}{2}\bigl| \gamma (\CMcal{D}^\prime \Omega )\phi \bigl| ,
\]
if $\phi \in \spinor _n$, in this case we used $\pi _{n+1}(e_j^{1,0}\cdot \nabla _{e_j}\pi _n)\pi _n=0$ and the fact that $\mathbf{i}+\gamma (\Omega )$ has absolute minimal eigenvalue $2$ on $\spinor _{n+1}^\perp $
The same method applied to $\pi _{n+1}\gamma (\Omega )=-\mathbf{i}\pi _{n+1}$ and $\gamma (\Omega )\pi _{n+1}=-\mathbf{i}\pi _{n+1}$ yields
\[
\Biggl| \sum _{j=1}^{2m}e_j^{0,1}(\nabla _{e_j}\pi _{n+1})\Biggl|\leq \frac{1}{2}|\CMcal{D}^{\prime \prime }\Omega |.
\]
Therefore, we obtain
\[
|\delta \frak{T}|\leq \frac{1}{2}\left( |\mathrm{d}^*\Omega |+|\CMcal{D}^\prime \Omega |+|\CMcal{D}^{\prime \prime }\Omega |\right) 
\]
which gives the claim $\widehat{\frak{R}}\geq 0$. 
\end{proof}
\section{Proof of the theorem}
\begin{lem}
Suppose $(M,g)$ is a complete spin manifold of real dimension $2m$. If the scalar curvature is uniformly bounded with $\mathrm{scal}\geq -4m(m+1)$, the Dirac operator
\[
\widetilde{\dirac}=\dirac -\mathbf{i}(m+1):W^{1,2}(M,\spinor M)\to L^2(M,\spinor M)
\]
is an isomorphism of Hilbert spaces.
\end{lem}
\begin{proof}
(cf.~\cite{AnDa,He3,MinO}) Using the Lichnerowicz formula proves that the bilinear form $
B(\varphi ,\psi )=\int _M\left< \widetilde{\dirac}\varphi ,\widetilde{\dirac}\psi \right> 
$ is coercive and bounded on $W^{1,2}(M,\spinor M)$. The surjectivity of $\widetilde{\dirac}$ follows from the Riesz representation theorem and \cite[Thm.~2.8]{GrLa3}. 
\end{proof}

Let $(M,g,J)$ be an almost Hermitian spin manifold which is strongly asymptotically complex hyperbolic, where $E\subset M$ is supposed to be the Euclidean end of $M$. We consider the connection $\widehat{\nabla }=\nabla +\frak{T}$ on $\spinor M$ and the connection $\widehat{\nabla }^0=\nabla ^0+\frak{T}^0$ on $\spinor M_{|E}$, where $\nabla ^0$ is the Levi--Civita connection and $\frak{T}^0$ is the K\"ahler Killing structure for the complex hyperbolic metric on $E$. The bundle $\bundle{V}^0\subset \spinor M_{|E}$ is trivialized by spinors parallel w.r.t.~$\widehat{\nabla }^0$.

The gauge transformation $A$ extends to a bundle morphism $A:\spinor M_{|E}\rightarrow \spinor M_{|E}$ with (cf.~\cite{AnDa})
\[
\left| \overline{\nabla }\varphi -\nabla \varphi \right| \leq C\left| A^{-1}\right| \left| \nabla ^{0}A\right| \left| \varphi \right| \, ,\]
where $\nabla $ is the usual spin connection for $g$ and $\overline{\nabla }$
is a connection on $\spinor M_{|E}$ obtained from the connection $\overline{\nabla }$ on $TM_{|E}$ and given by $\overline{\nabla }Y=A(\nabla ^{0}(A^{-1}Y))$.

Let $\psi _{0}$ be a spinor on $E\subset M$ which is parallel with respect to $\widehat{\nabla }^0$. Set $\psi :=h(A\psi _{0})$ for some cut off function $h$, i.e.~$h=1$ at infinity, $h=0$ in $M-E$ and $\mathrm{supp}(\mathrm{d}h)$ compact. We compute
\[ \begin{split}
\widehat{\nabla }_{X}\psi  = & (Xh)A\psi _{0}+h(\nabla _{X}A\psi _{0}+\frak {T}_{X}(A\psi _{0}))\\
 = & (Xh)A\psi _{0}+h(\nabla _{X}-\overline{\nabla }_{X})A\psi _{0}-hA\frak {T}^{0}_{X}\psi _{0}+h\frak {T}_{X}A\psi _{0}
\end{split}\]
and thus, the asymptotic assumptions supply
\[
\widehat{\nabla }\psi \in L^2(M,T^*M\otimes \spinor M)\]
and
\begin{equation}
\label{41}
\left\langle \widehat{\nabla }_{\nu }\psi +\nu \cdot \widehat{\dirac }\psi ,\psi \right\rangle \in L^1(M)
\end{equation}
($|\psi _0|_0^2$ can be estimated by $ce^{2r}$ for some $c>0$). Using the above lemma gives a spinor $\xi \in W^{1,2}(M,\spinor M)$ with $\widetilde{\dirac}\xi =\widetilde{\dirac}\psi \in L^2$. In particular $\varphi :=\psi -\xi$ is $\widetilde{\dirac}$--harmonic and non--trivial ($\psi \notin L^2$). Moreover, the selfadjointness of the boundary operator $\widehat{\nabla }_\nu +\nu \cdot \widetilde{\dirac}$ together with (\ref{41}) implies as usual
\[
\liminf _{r\to \infty }\int\limits _{\partial M_r}\left< \widehat{\nabla }_\nu \varphi +\nu \cdot \widetilde{\dirac}\varphi ,\varphi \right> =0
\]
for a non--degenerate exhaustion $\{ M_r\}$ of $M$ (cf.~\cite{AnDa}). Since inequality (\ref{ineq}) gives $\widehat{\frak{R}}\geq 0$, we conclude from the integrated Bochner--Weitzenb\"ock formula:
\[
\int\limits _{\partial M_r}\left< \widehat{\nabla }_\nu \varphi +\nu \cdot \widetilde{\dirac}\varphi ,\varphi \right> \geq \int\limits _{M_r}\left| \widehat{\nabla }\varphi \right| ^2\geq 0,
\]
that $\varphi $ is parallel w.r.t.~$\widehat{\nabla }$. Since $0=\widehat{\dirac}\varphi =\dirac \varphi -\mathbf{i}(m+1)\mathrm{pr}_\bundle{V}\varphi $ and $0=\widetilde{\dirac}\varphi =\dirac \varphi -\mathbf{i}(m+1)\varphi $, we obtain that $\varphi $ is a section of $\bundle{V}=\spinor _n\oplus \spinor _{n+1}$. Furthermore, $\widehat{\nabla }^0$ is a flat connection of $\bundle{V}^0$, so $\bundle{V}$ is trivialized by spinors parallel w.r.t.~$\widehat{\nabla }$. In particular, $\nabla _X$ preserves sections of $\bundle{V}$. Since $\widehat{\nabla }$ is flat on $\bundle{V}$, $\widehat{R}=0$ implies
\[
0=R_{X,Y}^s+[\frak{T}_X,\frak{T}_Y]+(\nabla _X\frak{T})_Y-(\nabla _Y\frak{T})_X
\]
on $\bundle{V}$. A straightforward computations shows that $(\nabla _X\frak{T})_Y$ and $(\nabla _Y\frak{T})_X$ are Hermitian on $\bundle{V}$ for all $X,Y$ (use the fact $(\frak{T}_X)^*=\frak{T}_X$), but $R_{X,Y}^s$ as well as $[\frak{T}_X,\frak{T}_Y]$ are skew--Hermitian on $\bundle{V}$ which leads to
\begin{equation}
\label{eq102}
0=R_{X,Y}^s+[\frak{T}_X,\frak{T}_Y].
\end{equation}
From the fact (cf.~\cite{BFGK})
\[ 
\gamma (\mathrm{Ric}(X)) =2\sum _{i}e_i\cdot R^s_{e_i,X}
\]
and equation (\ref{eq102}), we conclude $\mathrm{Ric}(X)=-2(m+1)X$ (cf.~\cite{Kir1}), i.e.~$g$ is Einstein of scalar curvature $-4m(m+1)$. Inequality (\ref{ineq}) yields $\mathrm{d}^*\Omega =0$ as well as $\CMcal{D}^\prime \Omega =0$ and $\CMcal{D}^{\prime \prime }\Omega =0$. In particular, $\CMcal{D}^\prime +\CMcal{D}^{\prime \prime }=\mathrm{d}+\mathrm{d}^*$ supplies $\mathrm{d}\Omega =0$. Therefore, if $J$ is integrable, $(g,J)$ must be K\"ahler (cf.~\cite[p. 148]{KoNo2}) and we could use the result by Herzlich to get the claim. However, we did not assume $J$ to be integrable and in order to prove the general case, we compute the Riemannian curvature of $(M,g,J)$. We have
\[ \begin{split}
[\frak{T}_X,\frak{T}_Y]=&(Y^{1,0}\cdot X^{0,1}\cdot -X^{1,0}\cdot Y^{0,1}\cdot )\pi _{n+1}+\\
&+(Y^{0,1}\cdot X^{1,0}\cdot -X^{0,1}\cdot Y^{1,0}\cdot )\pi _n\\
=&-\frac{1}{2}\gamma (X\wedge Y+JX\wedge JY)(\pi _n+\pi _{n+1})+\\
&+\mathbf{i}\Omega (X,Y)\pi _{n+1}-\mathbf{i}\Omega (X,Y)\pi _n
\end{split}\]
as well as $R^s_{X,Y}=\frac{1}{2}\gamma (\CMcal{R}(X\wedge Y))$, where $\CMcal{R}$ is the Riemannian curvature considered as endomorphism on $\Lambda ^2M$. Thus, we obtain
\[
\gamma \bigl(\CMcal{R}(X\wedge Y)-X\wedge Y -JX\wedge JY -2\Omega (X,Y)\Omega \bigl)\varphi =0
\]
for all $\varphi \in \Gamma (\CMcal{V})$ from (\ref{eq102}) and $\gamma (\Omega )\pi _n=\mathbf{i}\pi _n$, $\gamma (\Omega )\pi _{n+1}=-\mathbf{i}\pi _{n+1}$. In particular, the following lemma shows
\begin{equation}
\label{curvature}
\mathrm{pr}_{\Lambda ^{1,1}M}\circ \CMcal{R}(X\wedge Y)=X\wedge Y +JX\wedge JY +2\Omega (X,Y)\Omega .
\end{equation}
\begin{lem}
Suppose $(V,q)$ is a vector space of real dimension $2m$ with a qua\-dra\-tic form $q$ and a $q$--compatible complex structure $J$. Denote by $S=\oplus S_r$ the spinor space of $V$ where $S_r$ are induced from the action of the K\"ahler form $\Omega$. Choose $l:=\left[ \frac{m-1}{2}\right]$, then if $\eta \in \Lambda ^{1,1}V$ annihilates $S_l\oplus S_{l+1}$, i.e.
\[
\eta \cdot \psi =0
\]
for all $\psi \in S_l\oplus S_{l+1}$, $\eta $ has to vanish.
\end{lem}
\begin{proof}
Suppose that $m$ is even. The only $\Lambda ^{1,1}V$--forms which annihilate $S_{l+1}$ are multiples of $\Omega $ (cf.~\cite{He3}). But $\Omega $ acts as $2\mathbf{i}$ on $S_l$ which shows the claim if $m$ is even. Assume that $m$ is odd and $\eta \cdot \psi =0$ for all $\psi \in S_l\oplus S_{l+1}$. We consider the vector space $V\oplus \C ^2$ with its spinor space $S\widehat{\otimes }\C ^2$. Since Clifford multiplication with $\Lambda ^2V$ satisfies
\[
\omega \cdot (\psi \otimes \varphi )=(\omega \cdot \psi )\otimes \varphi ,
\]
and $(S\widehat{\otimes}\C ^2)_{l+1}$ is given by $S_l\otimes \C \oplus S_{l+1}\otimes \C $, we obtain $\eta \cdot \Psi =0$ for all $\Psi \in (S\widehat{\otimes}\C ^2)_{l+1}$. Thus, $\Lambda ^{1,1}V\subset \Lambda ^{1,1}(V\oplus \C ^2)$ together with the above case ($m$ even) implies that $\eta $ is a multiple of $\Omega $. But $\Omega $ acts as $\pm \mathbf{i}$ on $S_l$ respectively $S_{l+1}$ which shows the claim: $\eta =0$.
\end{proof}
Using equation (\ref{curvature}), the symmetry of the Riemannian curvature tensor and $\Omega \in \Gamma (\Lambda ^{1,1}M)$ lead to
\[
\left<\CMcal{R}(\Omega ),X\wedge Y\right>=\left< \CMcal{R}(X\wedge Y),\Omega \right> =2(m+1)\Omega (X,Y).
\]
Consider the Bochner--Weitzenb\"ock formula on $\Lambda ^2M$:
\[
\triangle =\mathrm{d}^*\mathrm{d}+\mathrm{d}\mathrm{d}^*=\nabla ^*\nabla +\frak{R} ,
\]
then $\frak{R}$ is given by $\mathrm{Ric}+2\CMcal{R}$ (cf.~\cite[Ap.~B]{Se2}), where $\mathrm{Ric}$ acts as derivation on $\Lambda ^2M$. We already know, that $g$ is Einstein, i.e.~$\mathrm{Ric}=-4(m+1)\mathrm{Id} _{\Lambda ^2M}$ supplies $\frak{R}(\Omega )=0$. Moreover, $\mathrm{d}\Omega =0$ and $\mathrm{d}^*\Omega =0$ imply that $\Omega $ is harmonic: $\triangle \Omega =0$, i.e.~we obtain $\nabla ^*\nabla \Omega =0$. Using the fact
\[
0=\triangle |\Omega |^2=\mathrm{d}^*\mathrm{d}|\Omega |^2=2\left< \nabla ^*\nabla \Omega ,\Omega \right> -2\left< \nabla \Omega ,\nabla \Omega \right> 
\]
we conclude that $(g,J)$ is K\"ahler. Thus, $\CMcal{R}:\Lambda ^2M\to \Lambda ^{1,1}M$ together with (\ref{curvature}) yield constant holomorphic sectional curvature $-4$ of $(M,g,J)$. Since the end of $M$ is diffeomorphic to $\R ^{2m}-\overline{B_R(0)}$, $M$ must be isometric to $\hyp{\C }^m$. 
\bibliographystyle{abbrv}
\bibliography{complex.bbl}

\end{document}